\documentclass[12pt]{article}\usepackage{amsfonts,amssymb}

\usepackage{amsmath,bbm}
\usepackage{theorem}          
\usepackage{latexsym}
\usepackage{float}                   % steuert Gleitobjekte wie Figuren,
\usepackage{times}
\usepackage{cite}

\newtheorem{tht}{Theorem}[section]

\newtheorem{thl}[tht]{Lemma}

\newtheorem{thc}[tht]{Corollary}

\newcommand{\mn}{\medskip}   

\newcommand{\sn}{\smallskip\noindent}

\newcommand{\cW}{{\mathcal{W}}} 
\newcommand{\cS}{{\mathcal{S}}}

\newcommand{\cN}{{\mathcal{N}}}
\newcommand{\cC}{{\mathcal{C}}}
\newcommand{\cX}{{\mathcal{X}}}
\newcommand{\cD}{{\mathcal{D}}} 
\newcommand{\cJ}{{\mathcal{J}}}
 
\newcommand{\Hh}{{\mathcal{H}}} 
\newcommand{\cG}{{\mathcal{G}}} 
\newcommand{\cK}{{\mathcal{K}}}

\newcommand{\dZ}{\mathbb{Z}}
\newcommand{\dR}{\mathbb{R}}
\newcommand{\dN}{\mathbb{N}}
\newcommand{\dC}{\mathbb{C}}

\newcommand{\gn}{\mathfrak{n}}
\newcommand{\gk}{{\mathfrak{k}}}

\newcommand{\gt}{{\mathfrak{t}}}

\newcommand{\im}{\mathrm{i}}

\begin{document}

\date{\small{Fakult\"at f\"ur Mathematik und 
Informatik\\ Universit\"at Leipzig, 
Augustusplatz 10, 04109 Leipzig, Germany\\ 
E-mail: schmuedgen@math.uni-leipzig.de }
}

\title{A Strict Positivstellensatz for the\\ Weyl Algebra}

\author{Konrad Schm\"udgen}

\maketitle

\renewcommand{\theenumi}{\roman{enumi}}
\begin{abstract}
Let $c$ be an element of the Weyl algebra $\cW (d)$ which is given 
by a strictly positive operator in the 
Schr\"odinger representation. It is shown that, under some conditions, 
there exist elements $b_1,{\dots}, b_d \in \cW (d)$ 
such that $\sum^d_{j=1} b_j c b^\ast_j$ is a finite sum of squares.

\end{abstract}
\section{Introduction}

In the last decade various versions and generalizations of the Archimedian Po\-si\-tiv\-stellensatz and of uniform denominator results have been obtained in semi-algebraic geometry (see the recent books \cite{PD},\cite{M1}). The proofs of these results are either purely algebraic \cite{R}, \cite{M2}, \cite{JP} or functional analytic \cite{S1}, \cite{PV}. The first proof of the Archimedean Positivstellensatz for compact semi-algebraic sets given in \cite{S1} was essentially based on methods from functional analysis. The purpose of this paper is to prove a strict Positivstellensatz for the Weyl algebra. Our approach uses again methods from operator theory and functional analysis. 

Let $d\in \dN$.
The Weyl algebra $\cW (d)$ (see e.g. \cite{D}) is the unital complex $\ast$-algebra with $2d$ hermitean generators $p_1,{\dots}, p_d, q_1,{\dots}, q_d$ and defining relations
\begin{align*}
&p_k q_k - q_k p_k = -\im \cdot \mathsf{1} ~~{\rm for}~~ k = 1, {\dots}, d,\\
&p_kp_l = p_l p_k, q_k q_l = q_l q_k, p_k q_l = q_l p_k~~{\rm for}~~ k, l = 1, {\dots}, d, k \ne l,
\end{align*}
where $\im$ denotes the complex unit and $\mathsf{1}$ is the unit element of $\cW (d)$. The Weyl algebra $\cW (d)$ has a distinguished faithful $\ast$-representation, the Schr\"odinger representation $\pi_0$. It acts on the Schwartz space $\cS (\dR^d)$, considered as dense domain of the Hilbert space $L^2 (\dR^d)$, by
\begin{equation*}
( \pi_0 (p_k) \varphi ) (t) = -\im \frac{\partial \varphi}{\partial t_k} (t), (\pi_0 (q_k) \varphi(t) = t_k \varphi (t), \varphi \in \cS (\dR^d), k = 1, 
{\dots}, d.
\end{equation*}
Setting $a_k := 2^{-1/2} (q_k + \im p_k), a_{-k} := 2^{-1/2} (q_k - \im p_k)$, the Weyl algebra $\cW (d)$ is the unital $\ast$-algebra with generators $a_1, {\dots}, a_k, a_{-1}, {\dots}, a_{-k}$, defining relations
\begin{align*}
&a_k a_{-k} - a_{-k} a_k = \mathsf{1}~~{\rm for}~~ k = 1, {\dots}, d\\
&a_k a_l = a_l a_k~~{\rm for}~~ k,l = -d, {\dots}, -1, 1, {\dots}, d, k \ne -l,
\end{align*}
and involution given by $a^\ast_k = a_{-k}, k = 1, {\dots}, d$. We abbreviate 
$$
N_k := a^\ast_k a_k~{\rm and}~N := N_1 + {\cdots} + N_d = a^\ast_1 a_1 + {\cdots} + a^\ast_d a_d.
$$

The Weyl algebra $\cW (d)$ has a natural filtration $(B_0, B_1,{\dots})$, where $B_n$ is the linear span of $a^{k_1}_1 {\cdots} a^{k_d}_d a^{l_1}_{-1} {\cdots} a^{l_d}_{-d}$ such that $k_1 + {\cdots} + k_d + l_1 + {\cdots} + l_d \le n$ and $k_j, l_j \in \dN_0$. Here, as usual, $a^0_j := \mathsf{1}$. The corresponding graded algebra associated with this filtration is the polynomial algebra $\dC [z, \overline{z}]\equiv \dC [z_1, {\dots}, z_d, \overline{z_1}, {\dots}, \overline{z_d}]$ in $2n$ complex variables $z_1, {\dots}, z_d, \overline{z_1}, {\dots}, \overline{z_d}$, 
where $z_j$ and $\overline{z_j}$ correspond to $a_j$ and $a^\ast_j$, respectively. If $c \in W_d$ is an element of degree $n$, we write $c_n (z,\bar{z})$ for the polynomial in $\dC [z, \bar{z}]$ corresponding to the component of $c$ with degree $n$.

Throughout this paper, $\alpha$ is a fixed positive number which is not an integer. Let $\cN$ denote the set of all finite products of elements $N + (\alpha + n) \mathsf{1} $, where $n \in \dZ$. Further, we shall use the set $\sum^2 (\cW (d))$ of all finite sums of elements $x^\ast x$, where $x \in \cW (d)$, and the positive cone
$$
\cW (d)_+ = \{x \in \cW (d) : \langle \pi_0 (x) \varphi, \varphi\rangle \ge 0~{\rm for~all}~\varphi \in \cS (\dR^d)\}.$$

The main result of this paper is the following
\begin{tht}\label{SI-T1}
Let $c$ be a hermitean element of the Weyl algebra $\cW (d)$ of even degree $2m$ and let $c_{2m} (z, \bar{z})$ be the polynomial of $\dC [z_1, {\dots}, z_d, \overline{z_1}, {\dots}, \overline{z_d} ]$ associated with the $2m$-th component of $c$. Assume that
\begin{enumerate}
\item[(i)] There exists $\varepsilon > 0$ such that $c - \varepsilon \cdot \mathsf{1} \in \cW (d)_+$.
\item[(ii)] $c_{2m} (z, \bar{z}) \ne 0~{\rm for~all}~z \in \dC^d, z \ne 0$.
\end{enumerate}
If $m$ is even, then there exists an element $b \in \cN$ such that $b c b \in \sum^2 (\cW (d))$. If $m$ is odd, then there exists $b \in \cN$ such that  $\sum^d_{j=1} b a_j c a^\ast_j b \in \sum^2 (\cW (d))$.
\end{tht}

{\em Note added in proof.} It can be shown that in Theorem 1.1 condition
(i) implies condition  (ii). That is, assumption (ii) in Theorem 1.1
can be omitted.

This paper is organized as follows. The proof of Theorem \ref{SI-T1} will be completed in Section \ref{S5}. In Sections \ref{S2} -- \ref{S4} we develop some technical tools. They are needed in the proof of Theorem 1, but they are also of interest in itself. In Section \ref{S2} we introduce and study algebraically bounded $\ast$-algebras. In Section \ref{S3} we define an auxiliary algebraically bounded $\ast$-algebra $\cX$ associated with the representation $\pi_0$ of the Weyl algebra. In Section \ref{S4} we classify the representations of this auxiliary $\ast$-algebra. The form of these representations is essentially used in the proof of  Theorem \ref{SI-T1} in Section \ref{S5}. A simple example illustrating Theorem 1.1 is presented in Section \ref{S6}.

Let us fix a few definitions and notations. By a $\ast$-{\it representation} \cite{S2} of a unital $\ast$-algebra $\cX$ on a pre-Hilbert space $\cD$ with scalar product $\langle\cdot,\cdot\rangle$ we mean an algebra homomorphism $\pi$ of $\cX$ into the algebra $L (\cD)$ of linear operators mapping $\cD$ into $\cD$ such that $\pi (\mathsf{1} ) = I$ and $\langle \pi (x) \varphi, \psi \rangle = \langle \varphi, \pi (x^\ast) \psi \rangle$ for $x \in \cX$ and $\varphi, \psi \in \cD$. Here $\mathsf{1}$  is the unit element of $\cX$ and $I$ is the identity map of $\cD$. The closure of an operator $y$ is denoted by $\bar{y}$. For a self-adjoint operator $y$, we  denote by $\sigma (y)$ the spectrum of $y$ and by $E_y (\cJ)$ the spectral projection of $y$ associated with a Borel set $\cJ$.

\section{The algebraically bounded part of a $\ast$-algebra}\label{S2}
In this Section $\cX$ is an arbitrary complex $\ast$-algebra with unit element $\mathsf{1}$. Let $\cX_h = \{ x \in \cX : x^\ast = x\}$ be the hermitean part of $\cX$. Each element $x \in \cX$ can be written as $x = x_1 + \im x_2$, where $x_1 \equiv Re\  x := \frac{1}{2} (x + x^\ast) \in \cX_h$ and $x_2 \equiv  {\rm Im}\  x := \frac{1}{2} \im (x^\ast - x) \in \cX_h$. Suppose that $\cX$ is an $m$-{\it admissible wedge} of $\cX$ in the sense of \cite{S2}, p.22, that is, $\cC$ is a subset of $\cX_h$ such that $\mathsf{1} \in \cC$, $x + y \in \cC, \lambda x \in \cC$ and $z^\ast x z \in \cC$ for all $x,y \in \cC, \lambda \ge 0$, and $z \in \cX$. Let $\succeq$ denote the ordering of the real vector space $\cX_h$ defined by $x \succeq y$ if only if $x - y \in \cC$. 

Let $\cX_b (\cC)$ be the set of all elements $x \in \cX$ for which there exists a positive number $\lambda$ such that
$$\lambda \cdot \mathsf{1} \succeq \pm Re\  x ~{\rm and}~ \lambda \cdot \mathsf{1} \succeq \pm \ {\rm{Im}}\  x.$$
\begin{thl}\label{S2L2}
\begin{enumerate}
\item[(i)] If $x, y \in \cX_b (\cC)$, then $x y \in \cX_b (\cC)$.
\item[(ii)] For $ x \in \cX_b$, we have $x \in \cX_b (\cC)$ if and only if $xx^\ast \in \cX_b (\cC)$. 
\item[(iii)] Let $x, y \in \cX_h$. If $x \succeq 0$ and $x - y = x y$, then $x \succeq y \succeq 0$.
\end{enumerate}
\end{thl}
{\bf Proof.} (i): We write $ x = x_1 + \im x_2$ and $y = y_1 + \im y_2$, where $x_1, x_2, y_1, y_2 \in \cX_h$. Since $x \in \cX_b (\cC)$ and $y \in \cX_b (\cC)$, there are positive number $\lambda$ and $\mu$ such that $\lambda \cdot \mathsf{1} \succeq  \pm x_j$ and $\mu \cdot \mathsf{1} \succeq \pm y_j$ for $j = 1,2$. Then, $\lambda \cdot \mathsf{1} \mp x_j \in \cC$. Therefore, by the definition of an $m$-admissible wedge, for $z \in \cX$ and $\alpha \in \dC$ we have
\begin{align*}
& (\alpha \cdot \mathsf{1} + z)^\ast (\lambda \cdot \mathsf{1} - x_1) (\alpha  \cdot \mathsf{1} + z) + (\alpha \cdot \mathsf{1} - z)^\ast (\lambda \cdot \mathsf{1} + x_1) (\alpha \cdot \mathsf{1} - z)\\
+~& (\alpha \im \cdot \mathsf{1} + z)^\ast (\lambda \cdot \mathsf{1} - x_2) (\alpha \im \cdot \mathsf{1} + z) + (\alpha \im \cdot \mathsf{1} - z)^\ast (\lambda \cdot \mathsf{1} + x_2) (\alpha \im \cdot \mathsf{1} - z)\\
=~& 4 \lambda (z^\ast z + | \alpha |^2 \cdot \mathsf{1}) - 2 \alpha z^\ast x^\ast - 2 \bar{\alpha} x z\in \cC.
\end{align*}
and hence
\begin{equation}\label{zx}
2 \lambda ( z^\ast z + |\alpha |^2 \cdot \mathsf{1}) \succeq \alpha z^\ast x^\ast + \bar{\alpha} x z.
\end{equation}

Setting $z^\ast =x$ and $\alpha = 2\lambda$ in (\ref{zx}) we get $4 \lambda^2 \cdot \mathsf{1}  \succeq x x^\ast$. Likewise, replacing $x$ by $y^\ast$ and $\lambda$ by $\mu$ we obtain $4 \mu^2 \cdot \mathsf{1}  \succeq y^\ast y$. In particular, the preceding proves the only if part of assertion (ii).
Setting now $z = y$ and inserting the relation $y^\ast y \preceq 4 \mu^2 \cdot \mathsf{1}$ just proved into 
(\ref{zx}), it follows that
\begin{equation}\label{alphayx}
2 \lambda (4 \mu^2 + | \alpha |^2) \cdot \mathsf{1}  \succeq \alpha y^\ast x^\ast + \bar{\alpha} xy.
\end{equation}
Letting $\alpha = \pm 1$ and $\alpha = \mp \im$ in (\ref{alphayx}), we conclude that
\begin{align*}
& \lambda (4 \mu^2 + 1) \cdot \mathsf{1} \succeq  \pm \frac{1}{2} (y^\ast x^\ast + x y) = \pm R e \ x y\\
& \lambda (4 \mu^2 + 1) \cdot \mathsf{1} \succeq  \pm \frac{1}{2} \im (y^\ast x^\ast - x y) = \pm Im \  x y.
\end{align*}
By definition the latter means that $x y \in \cX_b (\cC)$.\\
(ii): The only if part is already proven. It remains to show that $x x^\ast \in \cX_b (\cC)$ implies that $x \in \cX_b$. Since $x x^\ast \in \cX_b (\cC)$, there is a $\lambda > 0$ such that $\lambda \cdot \mathsf{1} \succeq  x x^\ast$. From the fact that 
$$
(x - \alpha \cdot \mathsf{1}) (x - \alpha \cdot \mathsf{1})^\ast = x x^\ast - \alpha x^\ast - \bar{\alpha}  x + 
|\alpha|^2 \cdot \mathsf{1} \in \cC
$$
it follows that
\begin{equation}\label{lambdax}
(\lambda + |\alpha|^2) \cdot \mathsf{1} \succeq   x x^\ast + |\alpha|^2 \cdot \mathsf{1}  \succeq   \alpha x^\ast  + \bar{\alpha}  x .
\end{equation}
Setting $\alpha = \pm 1$ and $\alpha = \pm \im$ in (\ref{lambdax}), we conclude that $R e\  x$ and $Im\  x$ are in $\cX_b (\cC)$ and so $x \in \cX_b (\cC)$.\\
(iii): From the relations $x = y + xy$ and $x \succeq  0$ we obtain $x y = y^2 + y x y \succeq  0$ and hence $ x= y + x y \succeq  y$. Using once more the assumptions $x \succeq  0$ and $y = x ( \mathsf{1} - y)$ we get $y - y^2 = ( \mathsf{1} - y) x ( \mathsf{1}-y)  \succeq  0$. Thus, $y \succeq  y^2 \succeq  0$. \hfill $\Box$
\begin{thc}\label{S2C3}
$\cX_b (\cC)$ is a unital $\ast$-subalgebra of $\cX$.
\end{thc}
{\bf Proof.} From its definition it is obvious that $\cX_b (\cC)$ is a $\ast$-invariant linear subspace of $\cX$. By Lemma \ref{S2L2}(i), $\cX_b (\cC)$ is a subalgebra of $\cX$.\hfill $\Box$

\sn
By the definition of $\cX_b(\cC)$ the unit element $\mathsf{1}$ is an order unit of the real ordered vector space $(\cX_b (\cC)_h,\succeq).$ The corresponding order unit seminorm $\| \cdot \|_1$ is defined by 
$$
\| x \|_\mathsf{1} = {\rm inf}~ \{ \lambda > 0 : \lambda \cdot \mathsf{1} \succeq x \succeq - \lambda \cdot \mathsf{1}\},~ x \in \cX_b (\cC)_h.
$$
Recall that a point $x$ is called an {\it internal point} of a subset $M$ of a real vector space $E$ if for any $y \in E$ there exists $\varepsilon > 0$ such that $x + \delta y \in M$ when ever $|\delta | < \varepsilon, \delta \in \dR$. Let $\cC^0_b$ denote the set of internal point of the wedge $\cC_b := \cC \cap \cX (\cC)_h$ in the real vector space $\cX_b (\cC)_h$. Clearly, $\cC^0_b$ coincides with the set of order units of $\cC_b$ in the order vector space $(\cX_b (\cC)_h, \succeq)$. In particular, $\mathsf{1} \in \cC^0_b$.
\begin{thl}\label{S2L3}
Let $z$ be an element of $\cX_b(\cC)_h$ which is not in $\cC^0_b$. Then there exists a state $F$ on the $\ast$-algebra $\cX_b (\cC)$ such that $F(z) \le 0$ and $F (x) \ge 0$ for $x \in \cC_b$.
\end{thl}
{\bf Proof.} Since $\cC^0_b$  is not empty, by Eidelheit's separation theorem for convex sets (see \cite{K}, \S 17, (3) or \cite{J}, 0.2.4) there exists a $\dR$-linear functional $f$ on $\cX_b (\cC)_h$ such that $f\not\equiv 0$ and $f(z) \le 0 \le f(x)$ for $x \in \cC_b$. Since $\mathsf{1} \in \cC^0_b$ and $f\not\equiv 0$, we have $f (\mathsf{1}) > 0$. We extend $f (\mathsf{1})^{-1} f$ on $\cX_b (\cC)_h$ to a $\dC$-linear functional $F$ on $\cX_b (\cC)$.\hfill $\Box$

\sn   
{\bf Remark 1.} From \cite{J}, 3.7.3 resp. 1.8.3, it follows that the $\cC_b$-positive state $F$ on $\cX_b(\cC)$ can be chosen to be extremal (that is, if $G$ is another state on $\cX_b (\cC)$ such that $0 \le G (x) \le F(x)$ for all $x \in \cC_b$, then $G=F$).

\smallskip
We now specialize to the case when $\cC$ is the $m$-admissible wedge $\sum^2 (\cX)$ of all finite sums of squares $x^\ast x$, where $x \in \cX$. In this case the $\ast$-algebra $\cX_b (\cC)$ is denoted by $\cX_b$ and called the {\it algebraically bounded part} of the $\ast$-algebra $\cX$. We say the $\ast$-algebra $\cX$ is {\it algebraically bounded} if $\cX = \cX_b$. The usefulness of these notions stems from the following obvious fact: For any $\ast$-representation $\pi$ $\ast$-algebra $\cX_b$ on a pre-Hilbert space $\cD$, each element $x \in \cX_b$ is mapped into a bounded operator $\pi (x)$ on $\cD$ and $\| \pi (x)\| \le \| x\|_{\mathsf{1}}$ for $x \in (\cX_b)_h$. Moreover, if the $\ast$-algebra $\cX$ has a faithful Hilbert space $\ast$-representation, then $\| \cdot \|_{\mathsf{1}}$ is a norm and the unit $\mathsf{1}$ is a inner point of the cone $\sum^2 (\cX_b)$ in the normed space $((\cX_b)_h, \|\cdot \|_{\mathsf{1}})$. 

\mn
We illustrate the preceding by a simple example which has been used in [PV].

\sn 
{\bf Example.} Let $\cX$ be the unital $\ast$-algebra generated by the rational functions $x_{kl} := x_k x_l(\mathsf{1} + x^2_1 + {\dots} + x^2_d)^{-1}, k, l = 0,{\dots}, d$, on $\dR^d$, where $x_0 :=\mathsf{1}$. Since all generators $x_{kl}$ are hermitean and 
$
0\le x^2_{kl} \le \sum^d_{i,j=1} x^2_{ij} =\mathsf{1},
$
it follows that $x^2_{kl} \in \cX_b$ and so $x_{kl} \in \cX_b$ by Lemma 2.1(ii). Hence the $\ast$-algebra $\cX$ is algebraically bounded.

\section{An auxiliary algebraically bounded $\ast$-algebra}  \label{S3}
%
%\subsection{}\label{S3SU1}
In what follows we use another unitarily equivalent form of the representation $\pi_0$, the so-called 
Fock-Bargmann representation (see e.g. \cite[1.6]{F}). For notational simplicity we shall write $x$ instead of $\pi_0(x)$ for $x \in \cW (d)$ and $\alpha$ instead $\alpha {\cdot} \mathsf{1}$ for $\alpha \in \dC$ when no confusion occurs. The Fock-Bargmann realization of the representation $\pi_0$ acts on the orthonormal basis $\{ e_\gn; \gn \in \dN^d_0\}$ of the representation Hilbert space by
\begin{equation}\label{fbrep}
a_k  e_\gn = n^{1/2}_k  e_{\gn - 1_k}, a_{-k} e_\gn = (n_k + 1)^{1/2} e_{\gn + 1_k}
\end{equation}
for $k = 1, {\dots}, d$ and $\gn = (n_1, {\dots}, n_d) \in \dN^d_0$. Here $1_k \in \dN^d_0$ denotes the multi-index with $1$  at the $k$-th place and zero otherwise and we set $e_{\gn-\mathsf{1}_k} = 0$ if $n_k =0$. The corresponding domain $\cD_0$ of the representation consists of vectors $\varphi = \sum_{\gn \in \dN_0^d} 
\varphi_\gn e_\gn$ such that $\sum_\gn n^r_1 {\dots}n^r_d |\varphi_\gn|^2 < \infty$ for all $r \in \dN$. Put $| \gn | := n_1 + {\cdots} + n_d$ for $\gn = ( n_1,{\dots}, n_d) \in \dN^d_0$. Then the actions of the elements $N_k$ and $N$ of the Weyl algebra are given by
\begin{equation}\label{nrep}
N_k e_\gn = n_k e_\gn~~{\rm and}~~ N e_\gn = |\gn | e_\gn, \gn \in \dN^d_0.
\end{equation}

Set $a_0 := 1$. We define the following operators on the domain $\cD_0$ :
\begin{align*}
& x_{kl} = a_k a_l (N + \alpha)^{-1}~{\rm for}~ k{=}0,{\dots}, d, l{=}-d, {\dots}, d; k{=} -d,{\dots}, d, l{=}0 , {\dots}, d,\\
&x_{-l,-k} = (N + \alpha)^{-1} a_{-l} a_{-k}~{\rm for}~ k,l = 0,{\dots}, d,\\
&x_k = x_{k0} = a_k (N + \alpha)^{-1} ~{\rm and}~ y_{k0} = x_{-k,k} = N_k (N + \alpha)^{-1}~{\rm for}~ k = 1, {\dots}, d,\\
&y_n = (N + \alpha + n)^{-1}~{\rm for}~ n \in \dZ.
\end{align*}
Let $\cX$ be the unital $\ast$-algebra generated by the operators $x_{kl}, k, l = -d, {\dots}, d$, and $y_n, n \in \dN_0$. The operator $x_{kl}, y_n$ resp. the $\ast$-algebra $\cX$ can be considered as non-commutative analogs of the Veronese map used in \cite{PV}. For $k, l = - d,{\dots}, d$ and $j = 1, {\dots},d$, we have
\begin{equation}\label{xklstar} 
x^\ast_{kl} = x_{-l,-k},\  x_{kl} = x_{lk}~{\rm if}~ k + l \ne 0,\  x_{j,-j} - x_{-j,j} = y_0.
\end{equation}
Note that the operators $y_n, n \in \dZ$, and $y_{k0}, k = 1, {\dots}, d$, pairwise commute. Moreover, $x_{ij} x_{kl} = x_{kl} x_{ij}$ for $i,j,k,l = 1, {\dots}, d$. From (\ref{fbrep}) and (\ref{nrep}) it is clear that all operators $x_{kl}, y_n$ and so all elements of $\cX$ are bounded on $\cD_0$ and leave $\cD_0$ invariant.

In order to formulate some relations we introduce the abbreviations $\gt (i,j) = 2$ if $ i > 0, j > 0, \gt (i,j) = 1$ if $i = 0, j >0$ or $i >0, j = 0$, and $\gt (i,j) = 0$ otherwise. For the rest of the paper we need a number of commutation relations of the operators defined above. They are easily verified by using formulas (\ref{fbrep}) and (\ref{nrep}). We shall list these relations in a convenient form for the applications given below. Not all relations are used in full strength.

\begin{align}
\label{r1}
&y_k - y_n = (n {-} k) y_k y_n = (n {-} k) y_n y_k ~~{\rm for}~~ k, n \in \dZ.\\
\label{r2}
&y_{10} + {\cdots} + y_{d0} = 1 - \alpha y_0.\\
\label{r3}
&x^\ast_{kj} x_{kj} = y_{k0} (y_{j0} - \delta_{kj} y_0), x^\ast_{k,-l} x_{k,-l} = (y_{k0} + \delta_{k l} y_0)(y_{l_0} + y_0)\nonumber\\
&{\rm for}~ j = 0, {\dots}, d, k ,\  l = 1, {\dots}, d.\\
\label{r4}
&y_0 x_{k l} = (1 + ({\rm sign} (k) + {\rm sign} (l)) y_0) x_{kl} y_0~~{\rm for}~~ k, l = -d, {\dots}, d.\\
\label{r5}
&y_n x^\ast_k = x^\ast_k y_{n+1}, x_k y_n = y_{n+1} x_k,\\
\label{r6}
&x_l x^\ast_k = x^\ast_k (1- y_2) x_l + \delta_{kl} y^2_1,\\
\label{r7}
&x_k x^\ast_k = y_{k0} y_1 (1 - y_1) + y^2_1, y_{k0} x^\ast_k = x^\ast_k (y_{k0} (1- y_1) + y_1),\\
\label{r8}
&x_{kl} y_0 = x_k x_l (1 - y_0), x_{-k,-l} y_0 = x^\ast_k x^\ast_l (1 + y_0),\\
\label{r9}
&x_{k,-l} y_0 = x_{-l,k} y_0 + \delta_{kl} y^2_0 = x^\ast_l x_k + \delta_{kl} y^2_0,\\
&{\rm for}~ ~k,l = 1, {\dots}, d ~~{\rm and }~~ n\in \dZ.\nonumber \\
\label{r10}
&x_{ij} x_{kl} - x_{kl} x_{ij} \in y_0 \cX, x_{ij} x_{kl} - x_{il} x_{kj} \in y_0 \cX,\\
\label{r11}
&y_0 a_k a_l = (1 + \gt (k,l) y_0) x_{kl},\\
&{\rm for}~~ i, j, k, l = -d, {\dots}, d.\nonumber
\end{align}
Moreover, we have $y_0 \cX = \cX y_0$.

\begin{thl} \label{5}
The $\ast$-algebra $\cX$ is algebraically bounded, that is, $\cX = \cX_b$.
\end{thl}
{\bf Proof.} From (\ref{r2}) and (\ref{r3}) we obtain
$$
(1 -\alpha y_0) y_0 = \sum^d_{k=1} y_{k0} y_0 = \sum^d_{k=1} x^\ast_{k0} x_{k0} \succeq 0$$
and $$y_0 = \alpha y^2_0 + \sum^d_{k=1} x^\ast_{k0} x_{k0} \succeq 0 ~{\rm and}~ \alpha^{-1} - y_0 = \alpha (y_0 - \alpha^{-1})^2 + \sum^d_{k=1} x^\ast_{k0} x_{k0} \succeq 0.$$
Therefore, we have
\begin{equation}\label{ynullpos}
\alpha^{-1} \succeq y_0 \succeq 0.
\end{equation}
Since $y_n - y_{n+1} = y_n y_{n+1}$ by (\ref{r1}), it follows from Lemma \ref{S2L2}(iii)  by induction on $n$ that $\alpha^{-1} \succeq y_n \succeq 0$ and so $y_n \in \cX_b$ for all $n \in \dN_0$. Using (\ref{r2}) and (\ref{r3}) we get
\begin{equation}\label{yjnull}
(1 - \alpha y_0)^2 = \left( \sum^d_{k=1} y_{k_0}\right)^2 = \sum_{k \ne l} x^\ast_{kl} x_{kl} + \sum^d_{k=1} y^2_{k0} \succeq y^2_{j0}
\end{equation}
for $j = 1, {\dots}, d$. Since $y_0 \in \cX_b$, from (\ref{yjnull}) and Lemma \ref{S2L2}(ii) we derive that $y_{j0} \in \cX_b$ for $j= 1, {\dots}, d$. Using (\ref{r3}) and Lemma \ref{S2L2}, (i) and (ii), it follows from the latter that $x_{kj} \in \cX_b$ for $k = 1, {\dots}, d$ and $j = - d,{\dots}, d$. Since $x_{-k,-j} = x^\ast_{jk}$, all generators of the $\ast$-algebra $\cX$ are in $\cX_b$. By  Corollary \ref{S2C3} (i), $\cX = \cX_b$. \hfill $\Box$

\smallskip
For the proof of Theorem \ref{SI-T1} below we need the following Lemma.
\begin{thl} \label{S5-L8}  For $n \in \dN$ and $i_1,{\dots},i_{4n} \in \{ -d,{\dots}, d\}$ there exist polynomials $f_j (y_0) \in \dR [y_0], j = 1,{\dots}, 2n$, such that $f_j(0)=0$ and 
\begin{equation}\label{ynulla}
y^n_0 a_{i_1}{\dots}, a_{i_{2n}} = f_1 (y_0) x_{i_1i_2} f_2 (y_0) {\cdots} f_n (y_0) x_{i_{2n-1} i_{2n}},
\end{equation}
\begin{equation}\label{aynull}
a_{i_{2n+1}} {\cdots} a_{i_{4n}} y^n_0 = x_{i_{2n+1} i_{2n+2}} f_{n+1} (y_0) {\cdots} f_{2n} (y_0) x_{i_{4n-1} i_{4n}},
\end{equation}
\begin{equation}\label{yay}
y^n_0 a_{i_1} {\cdots} a_{i_{4n}} y^n_0 = f_1 (y_0) x_{i_1 i_2} f_2 (y_0) {\cdots} f_{2n} (y_0) x_{i_{4n-1} i_{4n}}.
\end{equation}
\end{thl}
{\bf Proof.} It suffices to prove (\ref{ynulla}). Equation (\ref{aynull}) follows from (\ref{ynulla}) by applying the adjoint operation and (\ref{yay}) is obtained by multiplying (\ref{ynulla}) and (\ref{aynull}). 

We prove (\ref{ynulla}) by induction on $n$. For $n = 1$, formula (\ref{r11}) gives (\ref{ynulla}). We assume that (\ref{ynulla}) is true for $n$ and compute
\begin{align*}
y^{n+1}_0 & a_{i_1} {\dots} a_{i_{2n}} a_{i_{2n+1}} a_{i_{2n+2}} = y_0 f_1 (y_0) x_{i_1 i_2} {\cdots} f_n (y_0) x_{i_{2n-1}, i_{2n}} a_{i_{2n+1}} a_{i_{2n+2}}\\
&= \tilde{f}_1 (y_0) x_{i_1 i_2} {\cdots} \tilde{f}_n (y_0) x_{i_{2n-1}, i_{2n}} y_0  a_{i_{2n+1}} a_{i_{2n+2}}\\
&= \tilde{f}_1 (y_0) x_{i_1 i_2} {\cdots} \tilde{f}_n (y_0) x_{i_{2n-1},i_{2n}} (1 + \gt ( i_{2n+1}, i_{2n+2}) y_0) x_{i_{2n+1}, i_{2n+2}},
\end{align*}
where $\tilde{f}_j (y_0) \in \dR [y_0]$ and $\tilde{f}_j (0) = 0$. Here the first equality holds by the induction hypothesis. The second equality follows from (\ref{r4}), while the third one is obtained by inserting (\ref{r11}).
\hfill$\Box$

\section{Representations of the auxiliary $\ast$-algebra}\label{S4}
Suppose $\pi$ is an arbitrary $\ast$-representation of the $\ast$-algebra $\cX$ on a dense domain of a Hilbert space $\Hh$. Since $\cX = \cX_b$ by Lemma \ref{5}, all operators $\pi (x), x\in \cX$, are bounded, so $\pi$ extends by continuity to a $\ast$-representation, denoted again by $\pi$, on the Hilbert space $\Hh$. The aim of this section is to describe the structure of this representation $\pi$.
To shorten the notation, we write simply $x$ instead of $\pi (x)$ for $x\in \cX$ if no confusion is possible. Moreover, we use the multi-index notation
$$
x^\gn := x^{n_1}_1 {\cdots} x^{n_d}_d ~{\rm for}~ \gn = (n_1, {\dots}, n_d) \in \dN^d_0.
$$

\sn
{\bf {4.1}}
Let $\Hh_\infty := \ker y_0$ and let $\Hh_1$ be the closed linear span of subspaces $\cK_0 := \ker (y_0 - \alpha^{-1})$ and $\cK_\gn := (x^\gn)^\ast \cK_0$ f\"ur $\gn \in \dN^d_0$. In this subsection we show that $\Hh_\infty$ and $\Hh_1$ are invariant subspaces for the representation $\pi$ such that $\Hh = \Hh_1 \oplus \Hh_\infty$.

From the relations $y_0 y_n = y_n y_0, n \in \dN_0$, and (\ref{r4}) it is clear that $\Hh_\infty = \ker y_0$ is an invariant subspace for the representations $\pi$.
Since $y_n y_0 = y_0 y_n$ and $y_n x^\ast_k = x^\ast_k y_{n +1}$ by (\ref{r5}), $\cK_0$ and $\cK_\gn$ and hence $\Hh_1$ are invariant under $y_n, n \in \dN_0$. The invariance of $\Hh_1$ under $x^\ast_k$ is trivial.

Since $\alpha^{-1} \succeq y_0 \succeq 0$ by (\ref{ynullpos}),  the self-adjoint operator $y_0$ satisfies the relation $\alpha^{-1} I \ge y_0 \ge 0$ in the Hilbert space ordering. Hence its spectrum $\sigma (y_0)$ is contained in the interval $[0,\alpha^{-1}]$. 

Let $\varphi \in \cK_0$. Using the relations $(1 + y_0) y_1 = y_0$ by (\ref{r1}) and $x_k y_0 = y_1 x_k$ by (\ref{r5}), we have
$$
(1 + y_0) x_k \varphi = (1 + y_0) x_k (\alpha y_0 \varphi) = \alpha (1 + y_0)y_1 x_k \varphi = \alpha y_0 x_k \varphi$$
and  so $y_0 x_k \varphi = (\alpha - 1)^{-1} x_k \varphi$. Since $\sigma (y_0) \subseteq [0, \alpha^{-1}]$, the latter implies that
\begin{equation}\label{xkery}
x_k \varphi = 0 ~\ {\rm for}~\  \varphi \in \cK_0 = \ker (y_0 - \alpha^{-1}), k = 1, {\dots}, d.
\end{equation}

The invariance of $\cK_\gn$ and so of $\Hh_1$ under $x_k, k = 1, {\dots},d$, follows easily by induction on $|\gn|$ using relations (\ref{xkery}) and (\ref{r6}) and the fact that $\cK_\gn$ is invariant under $y_2$.

We prove the invariance of $\Hh_1$ under $x_{kl}$. Let $\varphi \in \cK_0$. Using (\ref{r1}) and (\ref{r5}) we compute
\begin{align}\label{xklx}
&x_{kl} (x^\gn)^\ast \varphi = \alpha x_{kl} (x^\gn)^\ast y_0 \varphi = \alpha x_{kl} (x^\gn)^\ast y_{|\gn|}(1 + |\gn| y_0) \varphi \nonumber \\
&= \alpha x_{kl} y_0 (x^\gn)^\ast (1+ |\gn| \alpha^{-1})\varphi
\end{align}
for $k, l = -d, {\dots}, d$. Expressing $x_{kl} y_0$ by means of relations (\ref{r8}) and (\ref{r9}) and using the invariance of $\Hh_1$ under $x_j, x^\ast_j$ and $y_0$, the right hand side of (\ref{xklx}) is in $\Hh_1$. Thus, the subspace $\Hh_1$  is invariant under the generators of $\cX$ and so under all representation operators.

We show that $\Hh_\infty \bot \Hh_1$. Indeed, if $\eta \in \Hh_\infty = \ker y_0, \varphi \in \cK_0 = \ker (y_0  - \alpha^{-1})$ and $n \in \dN^d_0$, then by (\ref{r1}) and (\ref{r5}) we have 
\begin{align*}
&\langle \eta, (x^\gn)^\ast \varphi \rangle = \langle \eta, \alpha (x^\gn)^\ast y_{|\gn|} (1+ |\gn| y_0)\varphi\rangle\\
&=\langle \eta, (\alpha + |\gn|)y_0 (x^\gn)^\ast  \varphi \rangle = \langle y_0 \eta, (\alpha + |\gn|))x^\gn)^\ast \varphi \rangle =0.
\end{align*}

Finally, we prove that $\Hh = \Hh_1\oplus \Hh_\infty$. Clearly $\cG := \Hh \ominus (\Hh_1 \oplus \Hh_\infty)$ is an invariant closed subspace for the representation $\pi$. We have to prove that $\cG = \{0\}$. Assume to the contrary that  $\cG \ne \{0\}$. Let $Y_0, Y_1$ and $X_k$ denote the restriction to $\cG$ of the operators $y_0, y_1$ and $x_k$ on $\Hh$, respectively. Since $\cG \bot \ker y_0$ and $\cG \bot \ker (y_0 - \alpha^{-1})$, we have $\ker Y_0 = \{0\}$ and  $\ker (Y_0 - \alpha^{-1}) = \{0\}$. Because $\sigma (Y_0) \subseteq \sigma (y_0) \subseteq [0, \alpha^{-1}]$, we therefore have $\lambda_0 := \sup \sigma (Y_0) > 0$. Fix $k \in \{1, {\dots},d \}$. By (10), $X_k Y_0 = Y_1 X_k$. This in turn implies that $X_k f (Y_0) = f (Y_0) X_k$ for all $f \in L^\infty (\dR)$ and so
\begin{equation}\label{xey}
X_k E_{Y_0} (\cJ) = E_{Y_1} (\cJ) X_k
\end{equation}
for any Borel subset $\cJ$ of $\dR$. Since $Y_1 = Y_0 (I + Y_0)^{-1}$ by (\ref{r1}), it follows from the spectral mapping theorem that $\lambda_0 (1 + \lambda_0)^{-1} = \sigma (Y_1)$. Because $\ker (Y_0 - \alpha^{-1}) = \{ 0\}$, for any $\varepsilon > 0$ there exists $\lambda \in \sigma (Y_0)$ such that $|\lambda - \lambda_0| < \varepsilon$ and $\lambda < \alpha^{-1}$. Hence we can choose numbers $\lambda_1 \in \sigma (Y_0)$ and $\delta > 0$ such that
\begin{equation} \label{lambda}
\lambda_0 (1 + \lambda_0)^{-1} < \lambda_1 - \delta < \lambda_1 + \delta \le \lambda_0,~ \lambda_1 + \delta < \alpha^{-1}.
\end{equation}
Let $\cJ := (\lambda_1 - \delta, \lambda_1 + \delta)$. Since $\lambda_1 \in \sigma (Y_0)$ and 
$\lambda_1 - \delta > \sup \sigma (Y_1)$, we have $E_{Y_0} (\cJ) \ne 0$ and $E_{Y_1} (\cJ) = 0$, so that $X_k E_{Y_0} (\cJ) =0$ by (\ref{xey}). Therefore, by (\ref{r3}) and (\ref{r2}),
$$
0 = \sum^d_{k=1} X^\ast_k X_k E_{Y_0} (\cJ) = (1 - \alpha Y_0) Y_0 E_{Y_0} (\cJ).$$
Because $\inf \{| (1 - \alpha \lambda) \lambda|; \lambda \in \cJ\} > 0$ by (\ref{lambda}) and $E_{Y_0} (\cJ) \ne 0$ we have obtained a contradiction. Thus, $\cG = \{0\}$ and $\Hh = \Hh_1 \oplus \Hh_\infty$.

\sn
{\bf{4.2}}
In this subsection we show that the restriction $\pi_1$ of the representation $\pi$ to $\Hh_1$ is a direct sum  of representations which are unitarily equivalent to the identity representation of $\cX$. By the identity representation we mean the representation $\rho$ of $\cX$ on the Hilbert space $\Hh_0$ given by $\rho (x) = \bar{x}, x \in \cX$, where $\bar{x}$ is the continuous extension of the operator $x$ on the dense domain $\cD_0$ to $\Hh_0$.\\
We begin with two preliminary lemmas.
\begin{thl}\label{S4-L6}
(i) $x_k (x^\gn)^\ast = (x^\gn)^\ast  ((1 - y_2) {\cdots} (1 - y_{|\gn|+1}))^2 x_k$ for all $k = 1, {\dots}, d$ and $\gn \in \dN^d_0, \gn \ne 0$, such that $n_k =0$.\\
(ii) $x_k x^\ast_k x^{\ast r}_k = x^{\ast r}_k (y_{k_0} (1 - (r + 1) y_{r+1}) + (r + 1) y_{r + 1} ) y_{r + 1}$ for $k = 1, {\dots}, d$ and $r \in \dN_0$.
\end{thl}
{\bf Proof.} (i) is proved by induction on $|\gn|$. If $|\gn| = 1$, then the assertion holds by (\ref{r6}). Suppose that the assertion is valid for $\gn$. Let $j \in \{ 1, {\dots},d\}, j \ne k$ and $\gn^\prime := \gn + \mathsf{1}_j$. Using the induction hypothesis and relations (\ref{r6}) and (\ref{r5}) we obtain
\begin{align*}
x_k (x^{\gn^\prime})^\ast &= x_k (x^{\gn})^\ast x^\ast_j = (x^\gn)^\ast ((1 - y_2){\cdots} (1 - y_{|\gn|+1}))^2 x_k x_j^\ast\\
&= (x^\gn)^\ast ((1 - y_2) {\dots} (1 - y_{|\gn|+1}))^2 x^\ast_j (1 - y_2) x_k\\
&= (x^{\gn^\prime})^\ast (( 1 - y_2) {\cdots} ( 1 - y_{|\gn^\prime| + 1} ))^2 x_k.
\end{align*}
(ii) is proved by induction on $r$. For $r = 0$ the assertion is just the first formula of (\ref{r7}). Suppose that the assertion holds for $r$. Using the induction hypothesis and relations (\ref{r5}) and (\ref{r7}) we compute
\begin{align*}
x_k x^\ast_k x^{\ast(r+1)}_k &= x^{\ast r}_k (y_{k0} (1 - (r + 1)y_{r+1}) + (r + 1) y_{r+1}) y_{r+1} x^\ast_k\\
&= x^{\ast r}_k (y_{k0}x^{\ast}_k (1 - (r + 1) y_{r + 2}) + (r+1) x^\ast_k y_{r + 2})y_{r + 2}\\
&= x^{\ast (r+1)}_k ((y_{k0} (1 - y_1)+ y_1) (1 - (r + 1) y_{r + 2}) + (r+1) y_{r + 2})y_{r + 2}\\
&= x^{\ast (r+1)}_k (y_{k0}(1 - (r + 2) y_{r + 2}) + (r+2) y_{r + 2})y_{r + 2},
\end{align*}
where the last equality is derived from relation (\ref{r1}). \hfill $\Box$
\begin{thl} \label{S4-L7}
If $\eta, \varphi \in \cK_0$ and $\gk,\gn \in \dN^d_0, |\gk| + |\gn| >0$, then 
\begin{equation} \label{etaphi}
\langle (x^\gk)^\ast \eta, (x^\gn)^\ast \varphi\rangle = \frac{n_1! {\cdots} n_d!}{((1 + \alpha) {\cdots} (|n| + \alpha))^2} \delta_{\gk,\gn} \langle \eta,\varphi\rangle
\end{equation}
\end{thl}
{\bf Proof.} First we prove that $(x^\gk)^\ast \eta \bot (x^\gn)^\ast \varphi$ if  $\gk \ne \gn$. Assume without loss of generality that $k_j > n_j$. Set $k^\prime_l = k_l, n^\prime_l = n_l$ for $l \ne j$, $ k^\prime_j = n^\prime_j =0$, and $\gk^\prime = (k^\prime_1, {\dots}, k^\prime_d), \gn^\prime = (n^\prime_1, {\dots}, n^\prime_d)$. From Lemma \ref{S4-L6}(i) is follows by induction on $s$ that there exists a polynomial $f$ (depending on $s$ and $\gn^\prime$) such that
\begin{equation}\label{xfrel}
x^s_j (x^{\gn^\prime})^\ast = (x^{\gn^\prime})^\ast f (y_2,{\dots}, y_{|\gn^\prime|+s} ) x^s_j~~{\rm for}~~s\in \dN.
\end{equation}
Further, using the formulas (\ref{r7}) it is easily shown by induction on $r$ that there exists a polynomial $g$ (depending on $r$) such that
\begin{equation}\label{xgrel}
x^r_j x^{\ast r}_j = g (y_{j_0}, y_1,{\dots},y_r)~~{\rm for}~~r\in \dN.
\end{equation}
Since $x_j y_{j_0} = (y_{j_0} (1 - y_1) + y_1) x_j$ by (\ref{r7}) and $x_j y_n = y_{n+1} x_j$ by (\ref{r5}) we conclude from (\ref{xfrel}) and (\ref{xgrel}) that there is a polynomial $h$ such that
\begin{equation}\label{xhrel}
x^s_j x^r_j x^{\ast r}_j (x^{\gn^\prime})^\ast = (x^{\gn^\prime})^\ast h (y_{j_0}, y_1,{\dots}, y_{r+s})x^s_j.
\end{equation}
Setting $s= k_j - n_j, r = n_j$ and using the fact that $x_j \varphi =0$ by (\ref{xkery}), (\ref{xhrel}) implies that $x^\gk (x^\gn)^\ast \varphi =  x^{\gk^\prime} x^s_j x^r_j x^{\ast r}_j (x^\gn)^\prime \varphi =0$ and so $ \langle (x^{\gk})^\ast \eta, (x^\gn)^\ast \varphi \rangle =0$.

Next we prove (\ref{etaphi}) in the case $\gk = \gn$. It clearly suffices to show that
\begin{equation}\label{xstarphi}
x^\gn (x^\gn)^\ast \varphi = \frac{n_1! {\cdots} n_d!}{((1 + \alpha) {\cdots} (|\gn | + \alpha))^2}~\varphi~~{\rm for}~~ \gn \in \dN^d_0, n \ne 0.
\end{equation}
We prove (\ref{xstarphi}) by induction on $|\gn|$. First we note that $y_{r+1} \varphi = (r+1 + \alpha)^{-1} \varphi$ by (\ref{r1}) and $0 = x^\ast_k x_k \varphi = x^\ast_{k0} x_{k0} \varphi = y_{k0} y_0 \varphi = \alpha^{-1} y_{k0} \varphi$ by (\ref{xkery}) and (\ref{r3}), so that $y_{k0} \varphi = 0$. Inserting these facts into Lemma \ref{S4-L6}(ii) we get
\begin{equation}\label{xjphi}
x_j x^\ast_j x^{\ast n_j}_j \varphi = x^{\ast n_j} (n_j + 1) (n_j + 1 + \alpha)^{-2} \varphi~~{\rm for}~~ n_j \in \dN_0.
\end{equation}
Setting $n_j = 0$, (\ref{xjphi}) gives (\ref{xstarphi}) for $|\gn| = 1$. Suppose that (\ref{xstarphi}) holds for $\gn$. Let $j \in \{1, {\dots}, d\}$. We prove that (\ref{xstarphi}) is true for $\gn^\prime = \gn + 1_j$. Set $\tilde{\gn} = (n_1, {\dots}, n_{j-1}, 0, n_{j+1}, {\dots}, n_d)$. Then we compute
\begin{align*}
x^{\gn^\prime} (x^{\gn^\prime})^\ast \varphi &=  x^\gn x_j x^\ast_j x^{\ast n_j}_j (x^{\tilde{^\gn}})^\ast \varphi = x^\gn x_j (x^{\tilde{^\gn}})^\ast  x^\ast_j x^{\ast n_j}_j \varphi\\
&=x^\gn (x^{\tilde{^\gn}})^\ast ((1-y_2){\cdots} (1-y_{|\tilde{\gn}|+1}))^2 x_j x^\ast_j x^{\ast n_j}_j \varphi\\
&=x^\gn (x^{\tilde{\gn}})^\ast ((1-y_2){\cdots} (1-y_{|\tilde{\gn}|+1}))^2 x^{\ast n_j}_j (n_j+1)(n_j + 1 + \alpha)^{-2} \varphi\\
&=x^\gn (x^{\tilde{^\gn}})^\ast ((1-y_{2+n_j}){\cdots} (1-y_{|\tilde{\gn}|+1 + n_j}))^2 (n_j+1)(n_j + 1 + \alpha)^{-2} \varphi\\
&=x^\gn (x^{\gn})^\ast  (n_j + 1)(|\gn|+1 + \alpha)^{-2} \varphi,
\end{align*}
where we used Lemma \ref{S4-L6}(i), formula (\ref{xjphi}) and the fact that $(1- y_k) \varphi = (k - 1 + \alpha)\\
(k + \alpha)^{-1}\varphi$. Inserting the induction hypothesis we obtain (\ref{xstarphi}) for $\gn^\prime$. \hfill$\Box$

\sn
Put $c_\gn := (n_1!{\cdots} n_d!)^{-1/2} (1 + \alpha) {\cdots} (|\gn| + \alpha)$ for $n \in \dN^d_0$, $\gn \ne 0$, and  $c_0 := 1$. Let $\{ \varphi_i; i \in I\}$ be an orthonormal basis of $\cK_0$. Then, by formula (\ref{etaphi}) the set  $\{ e_{\gn,i} {:=} c_\gn (x^\gn)^\ast \varphi_i; \gn \in \dN^d_0, i \in I\}$ is an orthonormal basis of $\Hh_1$. From 
$$
\| x^\ast_k (x^\gn)^\ast \varphi\| = (n_k + 1)^{1/2} (|\gn| + 1 + \alpha)^{-1} \| (x^\gn)^\ast \varphi\|, \varphi \in \cK_0,
$$
by (\ref{etaphi}) we derive
$$
x^\ast_k e_{\gn,i} = (n_k + 1)^{1/2} (|\gn| + 1 + \alpha)^{-1} e_{\gn + \mathsf{1}_k,i}.
$$
Therefore, by (\ref{fbrep}) and (\ref{nrep}), the operator $x^\ast_k$ acts on the orthonormal set  $\{e_{\gn,i}; \gn \in \dN^d_0\}$ as on the orthonormal basis $\{e_{\gn}; \gn \in \dN^d_0\}$ for the identity representation of $\cX$. The same  is true for the adjoint operator $x_k$ of $x^\ast_k$ and hence for all operators $y_n$ and $x_{kl}$ by (\ref{r8}) and (\ref{r9}). That is, for each $i\in I$ the restriction of $\pi_1$ to the closed linear span of vectors $\{e_{\gn,i}; \gn \in \dN^d_0\}$  is unitarily equivalent to the identity representation of $\cX$. Consequently, $\pi_1$ is the direct sum of representations of $\cX$ which are unitarily equivalent to the identity representation.

\sn
{\bf{4.3}}
In this subsection we study the restriction $\pi_\infty$ of $\pi$ to the invariant subspace $\Hh_\infty = \ker y_0$. Since $\pi_\infty (y_0) = 0$ and $x^\ast_{k0} = x^\ast_{0k} = x_{-k,0} = x_{0,-k}$, we have 
\begin{equation}\label{yinf}
\pi_\infty (y_n) = 0, n \in \dN_0,~~{\rm and}~~ \pi_\infty (x_{k0}) = \pi_\infty (x_{0k})= 0, k = -d,{\dots}, d.
\end{equation}
by (\ref{r1}) and (\ref{r3}). From (\ref{r10}), (\ref{r3}) and (\ref{xklstar}) we conclude that $X_{kl} := \pi_\infty (x_{kl}), k, l = -d,{\dots}, d$, are pairwise commuting bounded normal operators on $\Hh_\infty$ satisfying $X_{kl}= X_{lk}$, $X^\ast_{kl} = X_{-l,-k}$ and 
\begin{equation}\label{xijkl}
X_{ij} X_{kl} = X_{kj} X_{il}~~{\rm for}~~ i,j,k,l = -d,{\dots}, d.
\end{equation}
Recall that $y_{j0} = x_{j,-j}$. Therefore, by (\ref{r2}), 
\begin{equation} \label{xsum}
X_{1,-1} + {\cdots} + X_{d,-d} = I.
\end{equation}
For $j=1, {\dots},d$, we obtain from (\ref{xijkl}) and (\ref{xsum})
\begin{equation}\label{xlpos}
\sum^d_{k=1} X^\ast_{k,-j} X_{k,-j} = \sum^d_{k=1} X_{k,-k} X_{j,-j} = X_{j,-j}.
\end{equation}

We now describe  the Gelfand spectrum of the operator family $\{ X_{k,l}; k,l = -d, {\dots},d\}$ or equivalently  the character space of the abelian $C^\ast$-algebra generated by these operators. Let $\chi$ be such a character. From (\ref{xsum}), there is $j \in \{ 1,{\dots},d\}$ such that $\chi (X_{j,-j}) \ne 0$. Take $z_j \in \dC$ such that $z^2_j = \chi (X_{jj})$. Since $\chi (X_{j,-j}) \ge 0$ by (\ref{xlpos}) and $z^2_j z^{-2}_j = \chi (X_{jj} X_{-j,-j}) = \chi (X_{j,-j})^2$ by (\ref{xijkl}), we have $z_j\overline{z_j} = \chi (X_{j,-j}).$ 
We define $z_k := \chi (X_{k,-j}) \chi (X_{j,-j})^{-1} z_j$ for $k \ne j$. Note that the latter relation is trivially true for $k = j$, so it holds for all $k = 1, {\dots}, d$. Using the preceding facts and (\ref{xijkl}) we compute
\begin{align*}
z_k z_l &= \chi (X_{k,-j} X_{l,-j}) \chi (X_{j,-j})^{-2} \chi (X_{j,j})\\
&= \chi (X_{kl}) \chi (X_{-j,-j} X_{j,j}) \chi (X_{j,-j})^{-2} = \chi (X_{kl}),\\
z_k \overline{z_l} &= \chi (X_{k,-j}) (\overline{\chi(X_{l,-j})} \chi (X_{j,-j})^{-2} \chi (X_{j,-j})\\
&= \chi (X_{k,-j}X_{-l,j}) \chi (X_{j,-j}) = \chi (X_{k,-l})
\end{align*}
for $k,l = 1,{\dots},d$. From the latter and (\ref{xsum}) we get
$$
\sum^d_{k=1} z_k \overline{z_k} = \chi \left( \sum^d_{k=1} X_{k,-k}\right) = \chi (I) = 1.
$$
Thus we have shown that for each character $\chi$ there is a point $z = (z_1,{\dots}, z_d)$ of the unit sphere $S^d$ of the Euclidean space $\dC^d$ such that
$$
\chi (X_{kl})=z_kz_l~~{\rm and}~~\chi (X_{k,-l}) = z_k \overline{z_l}~~{\rm for}~~ k,l = 1,{\dots}, d.
$$
From the Gelfand theory it follows that there exists a spectral measure  $E(\cdot)$ on the unit sphere $S^d$ of $\dC^d$ such that
\begin{equation}\label{piz}
\pi_\infty (x_{kl}) = \int_{S^d} z_k z_l d E(z,\bar{z}), \pi_\infty (x_{k,-l}) = \pi (x_{-l,k}) = \int_{S^d} 
z_k \bar{z}_l d E(z,\bar{z})
\end{equation}
for $k,l = 1, {\dots},d$. Combined with (\ref{yinf}), these formulas describe the representation $\pi_\infty$ on the generators of $\cX$ completely.

\section{Proof of Theorem 1.1}\label{S5}
We first prove the assertion of Theorem \ref{SI-T1} in the case when $m$ is even, say $m {=} 2n$. Then $c \in \cW (d)$ has degree $4n$. From formula (\ref{yay}) in Lemma \ref{S5-L8} it follows that $y^n_0 c y^n_0$ belongs to  the $\ast$-algebra $\cX$.

The crucial step of the proof is to show that $y^n_0 c y^n_0 \in \sum^2 (\cX)$. Assume
% to the contrary that $y^n_0 c y^n_0 \not\in \sum^2 (\cX)$. 
the contrary. We apply Lemma \ref{S2L3} to the wedge $\cC = \sum^2 (\cX)$. Since $\cX = \cX_b$ by Lemma \ref{5}, there exists a state $F$ on the $\ast$-algebra $\cX$ such that $F(y^n_0 c y^n_0) \le 0$. Let $\pi_F$ denote the representation of $\cX$ with cyclic vector $\varphi_F$ associated with $F$ by the GNS construction such that $F(x) = \langle \pi_F (x) \varphi_F, \varphi_F\rangle$ for $x \in \cX$. As shown in Section \ref{S4}, $\pi_F$ decomposes into a direct sum of representations which are unitarily equivalent to the identity representation of $\cX$ on $L^2 (\dR^d)$ and the representation $\pi_\infty$ on $\Hh_\infty$. Let $\varphi_i \in L^2 (\dR^d), i \in I$, and $\varphi_\infty \in \Hh_\infty$ be the components of the vector $\varphi_F$ in this decomposition. Then, wen have
\begin{equation}\label{Fx}
 F(x) = \sum_{i \in I} \langle \bar{x} \varphi_i, \varphi_i\rangle + \langle \pi_\infty (x) \varphi_\infty, \varphi_\infty\rangle,\  x \in \cX.
\end{equation}
By assumption (i), $\langle y^n_0 c y^n_0 \varphi,\varphi \rangle = \langle c y^n_0 \varphi, y^n_0 \varphi\rangle \ge\varepsilon \langle y^n_0 \varphi, y^n_0 \varphi \rangle$ for $\varphi \in \cD_0 = \cS (\dR^d)$ and hence
\begin{equation}\label{ycy}
\langle \overline{y^n_0 c y^n_0} \varphi,\varphi\rangle \ge \varepsilon \| \overline{y^n_0} \varphi\|^2 > 0~{\rm for}~ \varphi  \in L^2 (\dR^d), \varphi \ne 0.
\end{equation}
From Lemma \ref{S5-L8} and the fact that $\pi_\infty (f_j (y_0))= \pi_\infty (f_j (0))$ we obtain 
$$
\pi_\infty (y^n_0 a_{i_1} {\dots} a_{i_{4n}} y^n_0) = \pi_\infty (x_{i_1 i_2}) {\cdots} \pi_\infty (x_{i_{4n-1}, i_{4n}})
$$
for $i_1,{\dots},i_{4n} = -d, {\dots}, d$. If the degree of a monomial $a_{i_1} {\cdots} a_{i_{4n}}$ is less than $4n$, then at least one index $i_j$ is zero and so $\pi_\infty (y^n_0 a_{i_1} {\cdots} a_{i_{4n}} y^n_0) =0$ by (\ref{yinf}). Hence we have $\pi_\infty (y^n_0 c y^n_0) = \pi_\infty (y^n_0 c_{4n} y^n_0)$. Using (\ref{piz})  we derive
\begin{equation}\label{cint}
\langle \pi_\infty (y^n_0 c y^n_0) \varphi_\infty, \varphi_\infty \rangle = \int_{S^d} c_{4n} (z,\bar{z}) d \langle E (z,\bar{z}) \varphi_\infty,\varphi_\infty \rangle.
\end{equation}
By assumption (ii), $c_{4n} (z,\bar{z}) >0$ for $ z\in S^d$. Since $F(y^n_0 c y^n_n) \le 0$, it follows from (\ref{Fx}), (\ref{ycy}) and (\ref{cint}) that all vectors $\varphi_i, i \in I$, and $\varphi_\infty$ are zero. But then $F(\mathsf{1})=0$ by (\ref{Fx}), in contradiction to the fact that $F$ is a state. Thus, $y^n_0 c y^n_0 \in \sum^2 (\cX)$.

That $y^n_0 c y^n_0 \in \sum^2 (\cX)$ means that there exist elements $g_1, {\cdots}, g_s \in \cX$ such that $y^n_0 c y^n_0 = \sum^s_{l=1} g^\ast_l g_l$. Let $b \in \cN$. Multiplying the latter equation by $b(N + \alpha)^n$ from the left and from the right we obtain
\begin{equation}\label{bcbsum}
bcb = \sum^s_{l=1} (g_l (N + \alpha)^n b)^\ast (g_l (N + \alpha)^n b).
\end{equation}
Each element of $\cX$ is a linear combination of finite products of operators $a_j$ and  $a^\ast_j,j = 1,{\dots}, d$, and $ y_k = (N + \alpha + k)^{-1}, k \in \dN_0$. Therefore, it follows from the relations $a_j y_k = y_{k+1} a_j$ and $a^\ast_j y_k = y_{k-1} a^\ast_j$ that we can choose $b \in \cN$ such that all denumerators $(N + \alpha + k)^{-1}$ of elements $g_l$ cancel, so that $g_l (N + \alpha)^n b \in \cW (d)$. Then we have $b c b \in \sum^2 (\cW (d))$ by (\ref{bcbsum}), as required.

Next we treat the case when $m$ is odd, say $m = 2n-1$. Then $\tilde{c} := \sum^d_{j=1} a_j c a^\ast_j$ has degree $4n$. By assumption (i) on $c$, we have 
$$
\langle \tilde{c} \varphi,\varphi \rangle = \sum^d_{j=1} \langle c a^\ast_j \varphi,a^\ast_j \varphi\rangle \ge \sum^d_{j=1} \varepsilon \langle a^\ast_j \varphi,a^\ast_j \varphi\rangle = \varepsilon \langle (N + d) \varphi, \varphi  \rangle \ge \varepsilon \langle \varphi,\varphi \rangle
$$
for $\varphi \in \cS (\dR^d)$. Since  $\tilde{c}_{4n} (z,\bar{z}) = c_{2m} (z,\bar{z})$ on $S^d$, $\tilde{c}$ satisfies assumptions (i) and (ii) too, so the preceding applies to $\tilde{c}$. This completes the proof of Theorem \ref{SI-T1}.

\sn
{\bf{Remark 2.}} The above proof shows that for even $m = 2n$ the assertion of Theorem \ref{SI-T1} remains valid if assumption (i) is replaced by the weaker requirement that the continuous extension of the bounded operator $(N + \alpha)^{-n} c (N + \alpha)^{-n}$ on $\cS (\dR^d)$ to $L^2 (\dR^d)$ is positive and has trivial kernel. The latter is satisfied if there exists a bounded positive self-adjoint operator $x$ on $L^2 (\dR^d)$ with trivial kernel such that $\langle c \varphi,\varphi) \ge \langle x \varphi,\varphi\rangle$ for $\varphi \in \cS (\dR^d)$. The special case $x = \varepsilon \cdot I$ is assumption (i).

\section{An Example}\label{S6}
Suppose that $d=1$. Since the spectrum of the closure of the operators $\pi_0 (N)$ is $\dN_0$ by (\ref{nrep}), a polynomial $p(N)$ of $N$ is in $W(1)_+$ if and only if $p(n) \ge 0$ for all $n\in \dN_0$. As shown in \cite{FS}, the element $p(N)$ belongs to $\sum^2 (W(1))$ if and only if there are polynomials  $q_0,{\dots}, q_{k} \in \dC [N], k \in \dN_0$, such that
\begin{equation}\label{ps}
p(N) = q_0 (N)^\ast q_0 (N) + Nq_1 (N)^\ast q_1 (N) + {\cdots} + N(N-1) {\cdots} (N-k+1) q_{k}(N)^\ast q_{k}(N).
\end{equation}
For $\varepsilon \ge 0$, we set $c_\varepsilon:=(N{-}1)(N{-}2) + \varepsilon$. From the preceding facts it follows that $c_\varepsilon$ is in $W(d)_+$ for all $\varepsilon \ge 0$ and that $c_\varepsilon$ is not in $\sum^2 (W(1))$ if $0\le \varepsilon < \frac{1}{4}$. Clearly, $c_\varepsilon$ satisfies the assumptions of Theorem 1.1 for all $\varepsilon >0$. For arbitrary real $\alpha$ we have
\begin{align*}
(N+\alpha)c_\varepsilon(N+\alpha)&= \mbox{$\frac{1}{2}$} \alpha^2 (N{-}1)^2 (N{-}2)^2 + (1+ \mbox{$\frac{1}{2}$}\alpha^2)N(N{-}1)(N{-}2)(N{-}3)\\
&~~~~ + (2\alpha+3)N(N{-}1)(N{-2})+ \varepsilon (N+\alpha)^2.
\end{align*}
The latter expression has been found by A. Sch\"uler. If $\alpha^2 < 2$, then the right hand side of the preceding equation is of the form (\ref{ps}) and so $(N+\alpha)c_\varepsilon(N+\alpha)\in \sum^2 (W(1))$ as asserted by Theorem 1.1.
\mn

\end{document}